\documentclass[twoside,10pt]{article}
\usepackage{amsmath, amsthm, amscd, amsfonts, amssymb, graphicx, color}
\usepackage{graphicx}
\usepackage{url}
\begin{document}
\setcounter{page}{1}
\setlength{\unitlength}{10mm}
\newcommand{\f}{\frac}
\newtheorem{theorem}{Theorem}[section]
\newtheorem{lemma}[theorem]{Lemma}
\newtheorem{proposition}[theorem]{Proposition}
\newtheorem{corollary}[theorem]{Corollary}
\theoremstyle{definition}
\newtheorem{definition}[theorem]{Definition}
\newtheorem{example}[theorem]{Example}
\newtheorem{solution}[theorem]{Solution}
\newtheorem{xca}[theorem]{Exercise}
\newtheorem{remark}[theorem]{Remark}
\numberwithin{equation}{section}
\newcommand{\e}{{\epsilon}}

\newcommand{\sta}{\stackrel}

\title{Conformally related  Douglas metrics are Randers}
\author {Vladimir  S. Matveev\thanks{Institut f\"ur Mathematik,
Friedrich-Schiller-Universit\"at Jena,
07737 Jena, Germany,  {vladimir.matveev@uni-jena.de}},  Samaneh  Saberali\thanks{Urmia University, Iran, samanehsaberali@gmail.com}}
\date{}
\maketitle
\begin{abstract}
\noindent We show that two-dimensional  conformally  related Douglas metrics are Randers
  
\end{abstract}
\section{Introduction}  A Finsler metric $F$ is called {\it Douglas}, if there exists an affine connection $\Gamma=(\Gamma_{jk}^i)$ such that each geodesic of $F$, after some re-parameterisation, is a geodesic of $\Gamma$. We assume w.l.o.g. that $\Gamma$ is torsion free. Such Finsler metrics were considered by J. Douglas  in \cite{D1,D2} and were  named Douglas metrics  (or metrics  of Douglas type) in    \cite{R}. 

Though results of our paper are local, let us note that 
partition of unity argument shows that the existence of such a connection locally, in a neighborhood of any point, implies its existence  globally. 

 Prominent examples of Douglas metrics are Riemannian metrics (with $\Gamma$ being  the Levi-Civita connection),
  Berwald metrics (in this case as  $\Gamma$    we can take the associated connection)  and locally projectively flat metrics (in this case, in the local coordinates such that the geodesics are straight lines one can take  the flat connection  $\Gamma\equiv 0$).

In the present paper we study the following question: {\it  can two conformally related Finsler metrics  $F $ and   $e^{\sigma(x)} F$ 
be both Douglas}? We do not require that the connection $\Gamma$ is the same for both  metrics, in fact by  \cite{S} two  conformally equivalent   metrics can not have the same (unparameterized) geodesics unless the conformal coefficient is constant.

 Of course two conformally related Riemannian metrics are both Douglas. Another trivial example is as follows: 
let   $F$ be Douglas and $\sigma$ be a constant.  Then,  $e^\sigma F$ is also Douglas.

 Let us give a less trivial example:
 \begin{example} \label{ex} 
Consider the Randers metric $F= \alpha + \beta$, where $\alpha(x,y)= \sqrt{g_{ij} y^iy^j} $ for a  Riemannian metric $g$ 
 and  $\beta$ is an 1-form. Assume in addition that $\beta$ is closed, locally it is equivalent to the condition that $\beta= df$ for a function $f$ on the manifold.  Then, the metric $F$ is Douglas, since adding the closed 1-form $\beta$ does not change the geodesics, so the geodesics  of $F$ are (up  to a  re-parameterisation)  geodesics of the Levi-Civita connection of $g$. 
 
Next, for any function $\gamma$ of one variable,  the conformally related metric $\tilde F= e^{\gamma(f(x))} F=  e^{\gamma(f(x))}  \alpha  + e^{\gamma(f(x))}  \beta$  is also Douglas.
 Indeed, since the 1- form $ e^{\gamma(f(x))}  \beta$  is closed, 
the geodesics of $\tilde F$ are geodesics of  $ e^{\gamma(f(x)) } \alpha$, i.e.,   geodesics  of the Levi-Civita connection of the   Riemannian metric  $e^{2 \gamma(f(x))} g$. \end{example}

It is easy to see (see e.g. \cite[Theorem 3.1]{GM})  that in the class of   Randers metrics the above example is the only possibly example (of conformally related Douglas metrics with nonconstant conformal coefficient).
 Indeed, Douglas  metrics are geodesically-reversible, in the sense that for any geodesic its
  certain orientation-reversing  unparameterisation is also a geodesic.   Now, it is known (see e.g. \cite[Theorem 1]{Matveev2012})  that for  a 
  geodesically reversible Randers metric  $\alpha+ \beta$ the 1-form $\beta$ is  necessary  closed. Thus, if two conformally related Randers metrics $\alpha + \beta$ and $ e^{\gamma( x )}\alpha+ e^{\gamma(x)}\beta$ are both Douglas, then  both 1-forms $\beta$ and $\gamma(x) \beta$ are closed, which locally 
   implies that $\beta= df$ and $\gamma$ is a function of $f$ as we claimed.   
 
Our main result  is that in dimension two only this example is possible:

\begin{theorem}\label{A}
 Let $F$ be a Douglas $2$-dimensional metric such that  the 
conformally related $e^{\sigma(x)}F$ is also Douglas. Assume $d\sigma\ne 0$ at a point $p$. Then, in a neighborhood of $p$ the metrics  $F$ and $\sigma F$ are as in Example \ref{ex} above: $F$ is Randers, $F= \alpha + \beta$, the 1-form $\beta$ is the differential of a function $f$, and
$\sigma$  is a function of $f$. 
\end{theorem} 
 
 Statements similar to Theorem \ref{A} appeared in literature before, in most cases one considered special Finsler metrics though. 
In particular,   \cite{GM} proves the analogous statement for  $ (\alpha,\beta)-$ metrics in dimension $n\geq3.$
 Conformally related Douglas $(\alpha, \beta)$ metrics were also considered in   \cite{SD}.  The question when 
   conformally related 
 Kropina metrics are Douglas was studied in \cite[Theorem 9]{M},
  where it was shown that  if 
 \underline{ every} conformal transformation of a Kropina metric $\alpha^{2}/\beta$ is Douglas, then     $\beta \wedge  d\beta= 0$, and vise versa. Note that  as it follows from   \cite{CMMM}, if a Kropina metric is Douglas, then $\beta \wedge  d\beta= 0$.

  Related results are \cite[Theorem 5]{V2},  \cite[Theorem 3]{V1} and \cite[Theorem 5]{V1}, where it is  proved that a  nontrivially conformally 
   related Berwald metrics are Riemannian, and also \cite[Theorem 8.1]{MT}, where an analogous statement was proved for Minkowski metrics. 
  
  All objects in our paper are assumed to be sufficiently smooth; Finlser mertics are assumed to be strictly convex. 
  
  {\em Acknowledgment:} We thank Julius Lang for useful discussions. VM was supported by DFG. Most results were obtained during the visit of SS to    Jena;  SS thanks  
Iranian ministry of education for financial support and  Friedrich Schiller University of Jena for hospitality.

  \section{Proof of Theorem \ref{A}}
  \subsection{Necessary conditions on  $F_{|T_xM}$ implied by the assumption that $F$ and $e^{\sigma(x)}F$ are both Douglas.} \label{sec:2.1}
  Recall that (arc-length parameterised) 
   geodesics of a Finsler metric $F$ are solutions  of the differential equation 
  \begin{equation}\label{eq:geodesics}
   \ddot x^i + 2  G^i = 0.
   \end{equation}
 Here $G^i= G(x,\dot x)^i$ are the so-called {\it spray coefficients}, they are given 
 by \begin{equation}\label{32}
G^{i}=\dfrac{1}{4}g^{il}\left([F^{2}]_{x^{k}y^{l}}-[F^{2}]_{x^{l}}\right),
\end{equation} 
 where $g_{ij}(x,y):=\dfrac{1}{2}[F^{2}]_{y^{i}y^{j}}(x,y)$  and $(g^{ij}):=(g_{ij})^{-1}.$   The notation $[ F^{2} ]_{x^k y^j}$ and later $F_{y^i y^j}$ means  the partial derivative with respect to indicated variables.    
 The condition that geodesics of $F$ are geodesics of the  affine connection $\Gamma$ is therefore equivalent to the condition 
 \begin{equation}\label{eq:proj}
 2 G^i(x,y)=  \Gamma_{jk}^i y^k y^j + P(x,y) y^i, 
 \end{equation} 
 which should be fulfilled for some function $P$ and 
 for any $y\in T_xM$. 
 
 Next, by replacing $F$  by $\tilde F= e^{\sigma(x)} F$ in \eqref{32}, we obtain 
 the following known relation  (see e.g. \cite[Equation (9.8)]{Y})
between the spray  coefficients $G^{i}$ of $F$ and $\tilde{G^{i}}$ of $\tilde F$:
\begin{equation}\label{si}
\tilde{G}^{i}=G^{i}+\sigma_{0}y^{i}-\tfrac{F^{2}}{2}\sigma^{i},
\end{equation} 
where $\sigma^{i} = g^{i\ell}\sigma _{\ell},$ $\sigma_i=\tfrac{\partial }{\partial x^i} \sigma$ 
  and $\sigma_0= \sigma_i y^i$. In view of this,  the condition that the metric $\tilde F$  is Douglas means the existence of 
   a  function $\tilde P$ and of a   torsion-free affine connection $\tilde \Gamma$ such that  
   \begin{equation}\label{eq:proj2}
 2 G^i(x,y)=  \tilde \Gamma_{jk}^i y^k y^j +  {F^2}  \sigma^i  + \tilde P(x,y) y^i.  
 \end{equation}
Combining (\ref{eq:proj}) and (\ref{eq:proj2}), we obtain 
\begin{equation}\label{eq:T}
T_{jk}^i y^k y^j  =    {F^2}  \sigma^i  + \hat P(x,y) y^i, 
\end{equation} 
where $T= \Gamma - \tilde \Gamma$ and $\hat P= (\tilde P -P)$.  Let us now multiply the equation $\eqref{eq:T}$ by $g_{is}$: we obtain: 
$$
T_{jk}^i g_{is} y^k y^j =   {F^2}  \sigma_s + \hat P(x,y) y^i g_{is}. 
$$
 In view of  $g_{is}= \tfrac{1}{2}[F^2]_{y^i y^s}=  FF_{y^i y^s}+ F_{y^i}F_{y^s}$,  this equation is equivalent to 
 $$
T_{jk}^i F F_{y^i y^s} y^k y^j   +  T_{jk}^i F_{y^i}F_{y^s} y^k y^j =   {F^2}  \sigma_s + \hat P(x,y)  y^i  F_{y^i} F_{y^s} +  \hat P(x,y)  y^i  F F_{y^i y^s}. 
$$
 Because $F$ is 1-homogeneous, $y^i  F_{y^i y^s}=0$ and $y^iF_{y^i}= F$. 
 Using these relations and rearranging the terms, we obtain 
\begin{equation} \label{-1} 
T_{jk}^i F F_{y^i y^s} y^k y^j   + F_{y^s} \left(T_{jk}^i F_{y^i} y^k y^j  - \hat P(x,y)   F \right)-F^2 \sigma_s=0. 
\end{equation}
We apply to the above equation the linear  operation 
$\xi_s \mapsto   \xi_s - \tfrac{1}{F} y^i \xi_i F_{y^s}$. 
In view of  $F F_{y^s}  - y^i F_{y^i} F_{y^s}=0$,   the middle term of the left hand side of  \eqref{-1} disappears and  the first term remains unchanged  in view of $y^iF_{y^i y^s}=0$.   After dividing the result by $F$, we obtain 
\begin{equation} \label{eq:main} 
T_{jk}^i F_{y^i y^s} y^k y^j =   F \sigma_s  - \sigma_0 F_{y^s}. 
\end{equation} 
 
 \begin{remark}
 We see that  \eqref{eq:main} does not contain derivatives with respect to $x$-variables. Thus, for a fixed $x\in M$, it is a system of partial differential equations on   $F$  restricted to   $T_xM$. 
 \end{remark} 
 
 \begin{remark}
 We see that the PDE-system  \eqref{eq:main} is linear in $F$. In fact, linearity was clear in advance due to the following geometrical argument: if the metrics $F_1$ and $F_2$ are Douglas with respect to the same connection $\Gamma$, and  the conformally related  with the same coefficient $e^\sigma$    metrics   $\tilde F_1 = e^\sigma F_1$ and $\tilde F_2= e^\sigma F_2$ 
 are also Douglas with respect to the same connection $\tilde \Gamma$, then   any  linear combination $\lambda_1 F_1 + \lambda_2 F_2$   
 is also Douglas (with respect to  $\Gamma$; we assume that $\lambda_1, \lambda_2 \in \mathbb{R}$ are such that $\lambda_1 F_1 + \lambda_2 F_2$  is a Finsler metric) and its conformally related  $(\lambda_1 F_1 + \lambda_2 F_2)e^\sigma$ is also Douglas (with respect to  $\tilde \Gamma$).
 \end{remark}

 \subsection{Proof of Theorem \ref{A}} 
 We assume that the function $\sigma(x)$ equals  $x^1$ (so in the  coordinates $d\sigma = (1,0)$), one can always achieve this by  a coordinate change. Next, we fix a point $x$ and work at the tangent space $T_xM= \mathbb{R}^2(y^1, y^2)$ at this point. Because of homogeneity, in the polar coordinates 
 $y^1= r \cos(\theta), y^2= r\sin(\theta)$ the function $F(y_1,y_2)$ is given by $r f(\theta)$.  
 The PDE-system \eqref{eq:main} reduces in this setting to one ODE on the function $f$; let us find this ODE.  By direct calculations we see that  the hessian (with respect to the coordinates $y^1, y^2$) of the function $F$ is given by  
\begin{equation} \label{convex}   \begin{pmatrix}
F_{y^{1}y^{1}} & F_{y^{1}y^{2}} \\
F_{y^{2}y^{1}} & F_{y^{2}y^{2}} 
\end{pmatrix} 
 = \frac{f(\theta)+f''(\theta)}{r}  \begin{pmatrix}
 \sin(\theta)^2 &   {-\cos(\theta) \sin(\theta)} \\
 {-\cos(\theta) \sin(\theta)} &   {\cos(\theta)^2} 
\end{pmatrix}.\end{equation}
Combining this with \eqref{eq:main}, we see that 
 \eqref{eq:main} is equivalent to the ODE 
 \begin{equation}\label{eq:ODE} (f''(\theta)+f(\theta))P(\theta)=   {f(\theta)\sin(\theta) +\cos(\theta) f'(\theta)}, 
\end{equation}
where  $$P(\theta)= K_0 \cos(\theta)^3 + K_1 \cos(\theta)^2 \sin(\theta) + K_2  \cos(\theta)  \sin(\theta)^2 + K_3    \sin(\theta)^3$$ and  where the constants $K_0,...,K_3$ are given by  
\begin{align*}
K_{0}=-T ^{2}_{11}\quad , K_{1}=T^{1}_{11}-2T^{2}_{12}, \quad K_{2}=2T^{1}_{12}- T_{22}^{2}, \quad  K_{3}=T^{1}_{22}.
\end{align*}

 Thus is a linear ODE of the 2nd order, so its solution space is at most two-dimensional. Actually, 
 locally, near the points where $P(\theta)\ne 0$,  it is precisely  two-dimensional, but not all local solutions can be extended to global solutions. Indeed,  since the variable $\theta$ ``lives'' on the circle,  we are only interested in $2\pi$-periodic solutions. Besides,  since the highest derivative of $f$ comes with the nontrivial  coefficient $P(\theta),$ a  solution can approach  
  infinity when   $\theta$ approaches $\theta_0$  such that  $P(\theta_0)= 0$.
   
  By direct calculations we see that  the function $\cos(\theta)$ is a solution.  \begin{remark} 
 Geometrically, the addition of the function $\cos(\theta)$ to a  solution $f$ 
  corresponds to the addition of the closed 1-form $dx^1$ to  $F = r f$; this operation does not change unparameterized  geodesics of the metrics $F$ and of the conformally related metric  $e^\sigma F + e^\sigma dx^1$, see the explanation in Example \ref{ex}. \end{remark}
 
 Our goal is to find all constants $K_0,...,K_3$ such that  there exists a  $2\pi$-periodic bounded solution $f$ of \eqref{eq:ODE}  such that
   $f$ is positive and  $f'' + f$  is positive (the last condition corresponds to the condition that $F$ is strictly convex, see e.g. \eqref{convex}).  We formulate the answer in the following lemma:

  \begin{lemma}\label{keylemma} 
  For the  constants $K_0,...,K_3$,  there exists a bounded solution $f(\theta)$ of ODE \eqref{eq:ODE}  such that it is $2\pi$-periodic and such that  $f$  and  $f'' + f$  are positive at all $\theta$ if and only if  \begin{align}\label{7}
 K_{0}=\frac{g_{12}g_{11}}{g_{11}g_{22}-g^{2}_{12}} ,\quad  K_{1}=1+\frac{3g^{2}_{12}}{g_{11}g_{22}-g^{2}_{12}},
 \end{align}
 \begin{align}\label{a}
  K_{2}=\frac{3g_{22}g_{12}}{g_{11}g_{22}-g^2_{12}}, \quad 
  K_3=\frac{g^2_{22}}{g_{11}g_{22}-g^2_{12}}
 \end{align}
 for a certain positively definite symmetric (constant) $2\times 2$-matrix  $g_{ij}$. In this case the general solution of  \eqref{eq:ODE} is given by 
 \begin{equation}\label{eq:gen}
 f(\theta)= \textrm{const}_1 \sqrt{\cos(\theta)^2 g_{11} + 2 \cos(\theta)\sin(\theta)g_{12} +  \sin(\theta)^2g_{22}}   + \textrm{const}_2 \cos(\theta). 
\end{equation}
  \end{lemma}

 Clearly, the solution  \eqref{eq:gen} corresponds to $F=  \textrm{const}_1 \alpha + \textrm{const}_2 \beta$, with $\alpha= \sqrt{g_{ij}y^iy^j}$ and $\beta= d\sigma= dx^1$; so Lemma \ref{keylemma}  together with explanation after Example \ref{ex}   imply Theorem \ref{A}.

 {\bf Proof of Lemma \ref{keylemma}.} The  direction ``$\Longleftarrow$'' (that for $K_0,...,K_3$ given by \eqref{a} the function \eqref{eq:gen} with $\textrm{const}_1>0$ and   $\textrm{const}_2=0$ is a $2\pi$-periodic  solution of \eqref{eq:ODE} satisfying the  assumptions in Lemma \ref{keylemma})  
 is clear geometrically and can be checked by calculations. Let us prove Lemma in the other (difficult) direction: we need to show  that the existence of such a solution $f$   implies that $K_0,...,K_3$ are as in \eqref{a}. We first replace the solution $f$ by  the function 
 $f_s(\theta):= f(\theta) + f(\theta + \pi)$. Since the function $P(\theta)$ satisfies $P(\theta + \pi)= -P(\theta)$, and also the functions $\cos(\theta)$ and  $\sin(\theta)$ satisfy $\cos(\theta + \pi)= -\cos(\theta)$, $\sin(\theta + \pi)= -\sin(\theta)$, all coefficients in the equation 
  \eqref{eq:ODE} change the sign after the addition of $\pi$ to the coordinate, so the function $ f(\theta + \pi)$, and therefore the function $f_s(\theta)$ is also a solution.  If $f$ is positive, $f_s$ is positive; if $f'' +f$ is positive, $f_s'' + f_s$ is positive. If $f$ is $2\pi$-periodic, $f_s$ is $\pi$-periodic, since $f_s(\theta + \pi)= f(\theta +\pi ) + f(\theta + 2 \pi)=f_s(\theta)$.

  Without  loss of generality we may and will  think that $f= f_s$, i.e., in addition to the above assumptions we also think that $f$ is $\pi$-periodic.

  \begin{remark} The operation $f \longrightarrow f_s$ corresponds geometrically to the ``symmetrisation'' $F(x,y) \longrightarrow F(x,y) + F(x,-y)$, this operation is compatible with  the conformal change of the metric and with the property of the metric to be Douglas with respect to  a connection $\Gamma$. \end{remark}

Observe that $(f'\cos(\theta) +f\sin(\theta))^{'}=\cos(\theta)(f+f'')$.  Denoting 
$(f'\cos(\theta)+f\sin (\theta))$ by $H$, 
  we obtain  $\frac{H'(t)}{\cos(\theta)} = \frac{H}{P(\theta)}$. 
  
   Since $ \frac{H'(\theta)}{\cos (\theta)}=f''(\theta)+f(\theta) $ and $ f''+f>0, $ we see that 
    \begin{equation} \label{def:G} G(\theta):= \frac{H'(\theta)}{\cos (\theta)}=\frac{H(\theta)}{P(\theta)}\end{equation} 
     is a smooth positive $\pi$-periodic function.  
    Its derivative  satisfies 
     $$G'(\theta)= \frac{H'(\theta) P(\theta) - H(\theta) P'(\theta)}{P(\theta)^2}\stackrel{\eqref{def:G}}{=} 
      \frac{G\cos(\theta) -G(\theta)P'(\theta)}{P(\theta)}.$$
      Which implies \begin{equation}\label{eq:ln} 
      \left(\ln(G(\theta))\right)' = \frac{\cos(\theta) - P'(\theta)}{P}.
      \end{equation}
   By our assumptions  the function $\ln(G(\theta))$ is smooth and $\pi$-periodic, so the following two conditions are satisfied: 
\begin{itemize}\item[(A)] $
 \int _{-\pi/2}^{\pi/2}\frac{\cos (\theta) -P'(\theta)}{P(\theta)}d\theta=0
$
 and \item[(B)]  $\frac{\cos(\theta)- P'(\theta)}{P(\theta)}$ is bounded. \end{itemize} 
 Our next goal is to  see  that   the existence of a solution of \eqref{eq:ln}  satisfying (A,B) is a strong condition on  $K_0,...,K_3$, 
 in fact we show   that $K_0,...,K_3$  such that there  exists  a solution of \eqref{eq:ln}  satisfying (A,B) are as in \eqref{a}.

 First, by direct calculations we observe 
  \begin{equation} \label{eq:tan} \frac{\cos (\theta)-P'(\theta)}{P(\theta)}=\frac{1-K_{1}-2K_{2}\tan( \theta)-3K_{3 } \tan^2(\theta)}{(K_{0}+K_{1}\tan(\theta)+K_{2}\tan^2(\theta)+K_{3}\tan^3(\theta))\cos^{2}(\theta)}+3\tan(\theta).\end{equation}

We consider this function restricted to the interval $\left(-\tfrac{\pi}{2}, \tfrac{\pi}{2}\right)$. There, $\tan(\theta)$ runs over  all real 
values, and the condition (B) implies that the polynomial  
   $ K_{0}+K_{1} t+K_{2} t^{2}+K_{3}t^{3} $ has only one root and this root is also a root of the polynomial  
  $1-K_{1}-2K_{2} t-3K_{3 }  t^2$.   Thus, \begin{eqnarray*} K_{0}+K_{1}\tan(\theta) +K_{2}\tan ^{2}(\theta)  +K_{3}\tan ^{3}(\theta) &=& (C+D\tan (\theta)+\tan ^{2}(\theta))(B-\tan(\theta)) E, \\  1-K_{1}-2K_{2} \tan(\theta)-3K_{3 }  \tan ^{2}(\theta) &=& 3(A-\tan(\theta))(B-\tan(\theta))E \end{eqnarray*} (for some constants $A,B,C,D, E$). 
  Then, the condition (A) reads (we make substitution $t= \tan(\theta)$ in the integral and also use that the function $\tan(\theta)$ is odd so $\int_{-r}^{r} 3  \tan(\theta) d\theta=0$ for each $r\in \left(-\tfrac{\pi}{2},  \tfrac{\pi}{2}\right)$.) 
  
   \begin{equation}
 \int_{-\pi/2}^{\pi/2} \frac{(A-\tan(\theta))}{(C+D\tan(\theta)+\tan^{2}(\theta))\cos^{2}(\theta)}dt= \int_{-\infty}^{\infty} \frac{(A-t)}{(C+Dt+t^2)  }dt=0.
\end{equation} 

Clearly, the   above integral is zero if and only if $C+Dy+y^{2}=(N+(A-y)^{2})$ for some constant $N$.  
In addition,  in order $(N+(A-y)^{2})$ to be nonzero (which is necessary by the condition (B)), we have $N=C-A^{2}>0$ (which in particular implies $C>0$).  
Thus, \begin{eqnarray} \label{kk1} K_{0}+K_{1} t+K_{2}t ^{2}+K_{3}t^{3} &=& ((C-A^{2})+(A- t)^{2})(B- t)E \\  1-K_{1}-2K_{2}t-3K_{3 } t ^{2}&=&3(A- t)(B-t)E. \label{kk2}\end{eqnarray}
Analyzing these two equations, 
we obtain $A=B$  and 
$E=\frac{1}{A^{2}-C}.$

 Combining  the formulas for $E,N,B$ as functions of $A,C$ obtained above with 
 (\ref{kk1}, \ref{kk2}), we  see that  $K_{0}, K_{1}, K_{2}, K_{3}$ are given by the following formulas:   
 
 \begin{align}\label{K0}
 K_{0}&=\frac{AC}{A^{2}-C}\\\label{k0}
 K_{1}&=1-\frac{3A^{2}}{A^{2}-C}\\\label{k1}
 K_{2}&=\frac{3A}{A^{2}-C}\\\label{k2}
 K_{3}&=\frac{-1}{A^{2}-C}
 \end{align}
By direct calculations we see that  the set of quadruples  $(K_0,...,K_3)$   given by the formulas (\ref{K0}--\ref{k2}) coincides  with the set of  quadruples  $(K_0,...,K_3)$  given by the formulas \eqref{7}. Indeed, for each $A$ and $C$, if we substitute in \eqref{7}   \begin{equation} \label{matrix}  g_{ij}=
\begin{pmatrix} C & -A \\  -A& 1\end{pmatrix},  \end{equation} 
 we obtain (\ref{K0}--\ref{k2}). Note that the condition $C-A^2=N>0$ implies that \eqref{matrix} 
is positively definite. Thus, the set of ``admissible'' quadruples $(K_0,K_1, K_2, K_3)$ (such that  \eqref{eq:ln} has a solution   satisfying (A,B)) is precisely the set of quadruples $(K_0,K_1, K_2, K_3)$ obtained by a symmetric positive definite matrix $g_{ij}$ by (\ref{K0}--\ref{k2}). 
Lemma \ref{keylemma}  and   Theorem \ref{A} are proved.

 \section{Prolongation of \eqref{eq:main} in higher  dimensions and conclusion  }
 
 Our initial goal  was to describe all nontrivially conformally related  Finsler 
 metrics such that both  are Douglas. We achieved this goal in dimension 2;  Theorem \ref{A} gives a complete answer. In higher dimension we do not know whether examples other than constructed in Example  \ref{ex} exist. In this section we would like to explain a way to approach a general problem (or to tackle the dimension 3 case).  We start with the following Remark:

 \begin{remark} \label{rem:3}
Suppose a function   $F=F(y)$ on $\mathbb{R}^n$  satisfies the following conditions: it is  positive, 
1-homogeneous, strictly convex  and  there exists a constant tensor  $T^i_{jk}$, symmetric with respect to the lower indexes,   and  a nonzero  constant covector $\sigma_i$ such that   \eqref{eq:main} are fulfilled. Then, we can     build a pair of nontrivially 
  conformally related  Douglas metric: the first one is the Minkowski metric given by $F_M(x,y)= F(y)$, 
  and the conformal related one is  $e^{\sum_i \sigma_i x_i}  F_M(x,y)$    (the function $\sigma(x)= \sum_i \sigma_i x_i$  chosen such that 
   its  differential is the constant covector  $\sigma_i$). 
   The associated connection of the metric $F_M$ is the flat one $\Gamma\equiv 0$,
    and of the metric $\sigma(x) F_M$ is given by $\tilde \Gamma^i_{jk}= -T^i_{jk}.$  
 \end{remark}

  The proof of the   statement  formulated in the remark  is straightforward and follows from the calculations in \S \ref{sec:2.1}: one simply needs to reverse all the arguments.   The only place where reversing the arguments may require additional  comments (because all others are in fact algebraic manipulations) is the transition from  \eqref{eq:main} to \eqref{-1}. Comparing these two formulas we see that 
  they are equivalent if $\left(T_{jk}^i F_{y^i} y^k y^j  - \hat P(x,y)   F \right)  =  - F\sigma_0 $; one can achieve it by choosing the appropriate function  $ \hat P(x,y)$. 
   
     Note  that  the Minkowski metric $F_M$ is clearly Douglas with  $\Gamma\equiv 0$. Now, by the construction of the equations \eqref{eq:main} we see that the difference between the spray coefficients of  
  $F$ and of $e^\sigma F$ is equal, up to the addition of  the appropriate  term of the form  $\check P(x,y) y^i$ (the  addition  of this term  does not change the geodesics),  to $T^i_{jk}$ so the metric $\tilde F= e^\sigma F$ is also Douglas with $\tilde \Gamma^i_{jk}= -T^i_{jk}.$

  We see that the  difficulty of our problem  is in fact  located in one tangent space; the dependence of the metric on the position is not important at least for the existence statement. Note that many previous researches working in this topic used the curvature of the metrics  and in particular  involved the  dependence of the metric on the point in the calculations; as we explained this is not necessary.

 Let us now study the equations \eqref{eq:main}. The system is clearly overdetermined; let us calculate the first compatibility conditions.  It can be done explicitly, the answer is given in the following Theorem. 
  
  \begin{theorem} \label{thm:prolongation}  In dimension $n\ge 3$, equations   \eqref{eq:main}  are fulfilled (for a certain 1-homogeneous smooth function $F$), if and only if   the following system of equations is  fulfilled 
  \begin{equation}\label{eq:main1}
 F_{y^i y^s} T_{jk}^i  y^k  -  F_{y^i y^j} T_{sk}^i  y^k   =   F_{y^j} \sigma_s  - \sigma_j F_{y^s}. 
  \end{equation}
  \end{theorem} 
  {\bf Proof.} The direction ``$\Longleftarrow$'' is easy: if we contract \eqref{eq:main1} with $y^j$ we obtain \eqref{eq:main}. Let us proof the statement in the direction ``$\Longrightarrow$''. Assume \eqref{eq:main} are satisfied. We differentiate them with respect to $y^\ell$ to obtain:
  \begin{equation}\label{eq:prolonation} 
F_{y^i y^s y^\ell } T_{jk}^i  y^k y^j  + 2 F_{y^i y^s  } T_{j\ell }^i   y^j =   F_{y^\ell} \sigma_s  - \sigma_\ell F_{y^s}- \sigma_0 F_{y^sy^\ell}. 
  \end{equation}
  Interchanging the indexes $\ell$ and $s$ in \eqref{eq:prolonation}, we obtain: 
 $$  F_{y^i y^s y^\ell } T_{jk}^i  y^k y^j  + 2 F_{y^i y^\ell   } T_{js  }^i  y^j =   F_{y^s } \sigma_\ell   - \sigma_s  F_{y^\ell}- \sigma_0 F_{y^sy^\ell}. $$
 Subtracting this equation from  \eqref{eq:prolonation}, we obtain:
 $$2 F_{y^i y^s  } T_{j\ell }^i   y^j - 2F_{y^i y^\ell  } T_{js  }^i   y^j = 2 F_{y^\ell} \sigma_s - 2F_{y^s } \sigma_\ell,$$
 which is clearly equivalent to \eqref{eq:main1}. Theorem \ref{thm:prolongation} is proved.

 Note that the number of  equations  in \eqref{eq:main1}, together with the 
 equations 
 \begin{equation} \label{eq:hom} F_{y^iy^s} y^i= 0\end{equation} 
 corresponding to the  1-homogeneity of $F$ is  $\tfrac{n(n+1)}{2}$  and is precisely the number of second derivatives of $F$ with respect to $y$ variables. It is easy to see that 
  for generic $T$ and in a neighborhood of almost every point $y$ 
 one can solve the system  with respect to the second derivatives and therefore to bring  the system in the Cauchy-Frobenius form, that is, all 
 highest derivatives of the unknown function $F=F(y)$ are  given as functions of the lower derivatives and of the coordinates $y$.   Indeed, since the system is linear, it is sufficient to show this for one tensor $T^i_{jk}$ (because determinate of a matrix is an algebraic expression in its components), and it is easy to find $T^i_{jk}$ such that the system has only one solution; one of the simplest examples is: 
 $$
 T^i_{jk}=\left\{ \begin{array}{ccl}0 & \textrm{ if } &  i\ne 1 \\
 0 & \textrm{ if } &  j\ne k \\
 1 & \textrm{ if } &  i=1  \textrm{ and  } j =k
   \end{array}\right. . $$ (This example corresponds to   the flat metric $g_{ij}=\delta_{ij}$ and to the 1-form  $\sigma_i=(1,0,...,0)$).    
  Now, from the  general theory it follows that for  $T^i_{jk}$  such that the solution on the system  is unique  the restriction of the Finsler metric to the tangent space    $T_xM$    depends on finitely many parameters (which in our case are the values of the first $y$-derivatives at one point of $T_xM$).

  Combining Theorem \ref{thm:prolongation} with Remark \ref{rem:3}, we obtain: 
  
  \begin{theorem} Assume there  exists a  convex 1-homogeneous  strictly convex  function  $F:\mathbb{R}^n\to \mathbb{R}$ satisfying (\ref{eq:main1},\ref{eq:hom}) such that the level set $\{y\mid F(y)= 1\}$ is not an ellipsoid. Then, there exists a nonranders metric $F$  
  such that it is Douglas and a conformal coefficient such that the  conformally related is also Douglas. 
  \end{theorem}

 We do not know whether such functions exists. The system (\ref{eq:main1},\ref{eq:hom}) is a linear overdetermined 
 system of PDE and in theory there exists an algorithmic was to solve  it, but  we did not managed to go through algebraic difficulties in the case of general $T^i_{jk}$. Besides, it is not clear how to analyze whether the  solution is indeed strictly convex (note that in dimension two strict convexity corresponds to the condition $f'' + f>0$ and was essentially used in the proof), and also how to analyze the solutions near the points where the solutions of the system   (\ref{eq:main1},\ref{eq:hom})  is not unique (in dimension two the analog of such points are the points where the coefficient of $f''$ vanishes; they play important role in the proof).  But for explicitly given ``special metrics'' (i.e., $(\alpha,\beta)$ metrics), in order to understand whether one can construct nontrivially conformally equivalent Douglas metrics in their  class, one simply should put the form of the metric in \eqref{eq:main1}  and then analyze the obtained equations.

\end{document}